\documentclass[12pt]{amsart}
\usepackage{amssymb, amsfonts, amsmath, fullpage, epsfig}
\usepackage{comment}

\newcommand{\noi}{\noindent}

\parskip=\smallskipamount
\numberwithin{equation}{section}

\newtheorem{izr}{Theorem}[section]
\newtheorem{lm}[izr]{Lemma}
\newtheorem{pos}[izr]{Corollary}

\theoremstyle{definition}

\newtheorem{op}[izr]{Remark}

\newcommand{\Nn}{\mathbb{N}}

\newcommand{\Cc}{\mathbb{C}}

\newcommand{\Zz}{\mathbb{Z}}

\newcommand{\dok}{\noindent \it Proof. \rm }
\newcommand{\ra}{\rightarrow \:}

\newcommand{\ol}[1]{\overline{#1}}

\newcommand{\psh}{plurisubharmonic }
\newcommand{\nbhd}{neighbourhood }
\newcommand{\va}{\alpha}
\newcommand{\vb}{\beta}
\newcommand{\vph}{\varphi}

\newcommand{\dbar}{\overline{\partial}}

\newcommand{\sgn}{{\rm sign}}
\newcommand{\m}[1]{\mathcal #1}

\def\dim{\mathop{\rm dim}\nolimits}

\def\Tr{\mathop{\rm Tr}\nolimits}

\def\inic{\setcounter{izr}{0}}

\begin{document}

\inic \pagestyle{plain} \setcounter{section}{0} \setcounter{tocdepth}{1}

\title{$\dbar$-equation on $(p,q)$-forms on conic neighbourhoods of $1$-convex manifolds}
\author{Jasna Prezelj}
\address{Jasna Prezelj, Faculty of Mathematics and Physics, Department of Mathematics,
University of Ljubljana, Jadranska 21, SI-1000 Ljubljana,
Slovenia}
\address{Faculty of Mathematics, Natural Sciences and Information Technologies, University of Primorska, Glagolja\v ska 8, SI-6000 Koper, Slovenia}
\email{jasna.prezelj@fmf.uni-lj.si}
\thanks{The author was supported by research program P1-0291 and by research project J1-5432 at Slovenian Research Agency. Part of the paper was written while author was visiting NTNU, Trondheim, Norway, and she wishes to thank this institution for its hospitality.}

\bigskip\bigskip\rm
\begin{abstract} Let $X$ be a $1$-convex manifold with the exceptional set $S$, which is also a manifold, $Z$ a complex manifold, $Z \ra X$ a holomorphic submersion,  $a: X \ra Z$ a holomorphic section and $S \subset U \Subset X$ an open relatively compact $1$-convex set. We construct a metric on a vector bundle $E \ra Z$ restricted to  a neighbourhood $V$ of  $a(U),$ conic along $a(S)$  with at most polynomial poles at $a(S)$ and positive Nakano curvature tensor in bidegree $(p,q).$
\end{abstract}


\maketitle
\section{Introduction and the main theorem}

\begin{figure}[h!]
\begin{center}
 \epsfysize=45mm
 \epsfbox{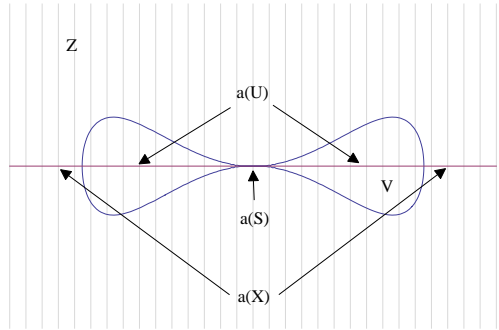}
 \caption{Conic neighbourhoods of $a(U \setminus S)$ in $Z$}
 \label{okolice}
\end{center}
  \end{figure}

\noi Let $\pi: Z \ra X$ be a submersion from a complex manifold $Z$ to  a $1$-convex manifold $X$ with the exceptional set $S,$ which is also a manifold.
The  motivation for this work was  to construct
solutions of the $\dbar$-equation with at most polynomial poles at $\pi^{-1}(S)$ in a  particular geometric situation (see Fig. \ref{okolice}), namely on a conic neighbourhood $V$ of sections of  $\pi: Z \ra X.$ The results of this paper may provide one step in the proof of such a claim.
It turned out that it is  possible to construct solutions of the $\dbar$-equation such that their $L^2$-norms on $V_{\delta}:= \lbrace z \in V, d(z,\pi^{-1}(S))> \delta\rbrace$  with respect to some ambient Hermitian metric $h_Z$ on $Z$ and $h_E$ on $E$ grow at most polynomially  as $\delta \ra 0.$   To this end weights that give rise to positive Nakano curvature and have polynomial behaviour on $\pi^{-1}(S)$ are constructed. Sup-norm estimates are given in the case $q = 1.$ The main theorem  of the present paper is the following:

\begin{izr}[Nakano positive curvature tensor in bidegree $(p,q)$] \label{main thm1}
  Let $Z$ be a $n$-dimensional complex manifold, $X$  a $1$-convex manifold, $S \subset X$ its exceptional set, which is also a manifold,
  $\pi : Z \ra X$ a holomorphic submersion with $r_0$-dimensional fibres, $\sigma : E \ra Z$ a holomorphic vector bundle and $a: X \ra Z$
  a holomorphic section. Let $\vph : X \ra [0,\infty)$ be  a plurisubharmonic exhaustion function, strictly plurisubharmonic on $X \setminus S$ and $\vph^{-1}(0)  = S.$ Let $U  = \vph^{-1}([0,c])$ for some $c > 0$ be a given holomorphically convex set and let $s \in \lbrace 1,\ldots,n \rbrace.$ Then there exist an open neighbourhood $V$ of $a(U\setminus S)$ conic along $a(S),$ a Nakano positive Hermitian metric $h$ on  $E|_V$ with at most polynomial poles on $\pi^{-1}(S)$ and such that the Chern curvature tensor $i\Theta(E \otimes \Lambda^s TZ)$ restricted to $V$ is Nakano positive and has at most polynomial poles and zeroes on $\pi^{-1}(S).$
\end{izr}

 For standard techniques for solving the $\dbar$-equation we  refer  to Demailly's book Complex analytic and algebraic geometry \cite{dem}. Our main tool is
 \begin{izr}[\rm Theorem VIII-4.5, \cite{dem}]\label{dem1} If $(W, \omega)$ is complete and $A_{E,\omega}> 0$ in bidegree $(p, q)$, then for any $\dbar$-closed form  $u\in L^2_{p,q} (W,E)$ with
$$
  \int_W \langle A^{-1}_{E,\omega} u, u \rangle dV < \infty
$$
there exists $v \in L^2_{p,q-1}(W,E)$ such
that $\dbar v = u$ and
$$
  \|v\|^2  \leq  \int_W \langle A^{-1}_{E,\omega} u, u \rangle dV.
$$
\end{izr}
It seems difficult to handle the commutator $A_{E,\omega}$ (see Sect. VII-7, \cite{dem} for computations) in the case of $(p,q)$-forms because it has terms of mixed signs for $p < n.$ Therefore
we view the $(p,q)$-form $u$ as a $(n,q)$-form $u_1$ with coefficients in  $E_1 = E \otimes \Lambda^{n-p}TZ $ by invoking the isomorphism $\Lambda^p T^*Z \simeq \Lambda^{n-p} TZ \otimes \Lambda^n T^*Z.$ If $u$ is closed so is $u_1.$  In addition,  in bidegree $(n,q),$
there is an analog of Theorem \ref{dem1} for a noncomplete K\"{a}hler metric provided that the manifold possesses a complete one (Theorem VIII-6.1, \cite{dem}). Moreover, the positivity of $A_{E_1,\omega}$ follows from the positivity of the Chern curvature tensor $i \Theta(E_1) = i\Theta (E) + i \Theta( \Lambda^{n-p}TZ ).$

\begin{figure}[h!]
\begin{center}
 \epsfysize=45mm
 \epsfbox{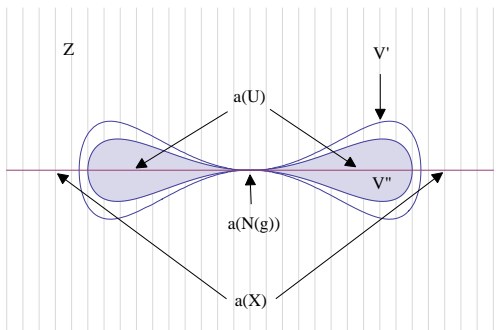}
 \caption{Conic neighbourhoods of $a(U \setminus N(g))$ in $Z$}
 \label{okoliceb}
\end{center}
  \end{figure}

To use both theorems we first need a Nakano positive Hermitian metric on $E|_V$ and a K\"{a}hler metric on $V,$ both with polynomial zeroes or poles on  $\pi^{-1}(S)$. The K\"{a}hler metric $\omega = i\partial \dbar \Phi$ constructed in Sect. 2 and the Nakano positive Hermitian metric   given by Theorem 1.1 in \cite{pre2} have the desired properties.  The space $(V,\omega)$ is  not complete but if we take a smaller neighbourhood conic  along $a(N(g)),$ $N(g):=g^{-1}(0),$ for some holomorphic function $g : X \ra \Cc$ with  $g(S) = 0,$ it contains a conic Stein \nbhd $V'$ which is complete K\"{a}hler (Fig. \ref{okoliceb}). With the above notation  we have

\begin{pos}[$\dbar$-equation in bidegree $(p,q)$] \label{pqforme}  Let $g : X\ra \Cc$ be holomorphic with $N(g) = g^{-1}(0) \supset S.$ Let $V'$ be an open Stein neighbourhood of $a(U \setminus N(g)),$ conic along $a(N(g)),$  let $u$ be a closed $(p,q)$-form on $V',$  $u_1$ the corresponding $(n,q)$-form with coefficients in $E_1 = E \otimes \Lambda^{n-p}TZ $ and  let $h$ be the metric on $E$ from  Theorem \ref{main thm1} with $s = n-p$ and  $h_1$ the induced metric on $E_1.$ Denote by $A = A_{E \otimes \Lambda^{n-p}TZ, \omega}$  the commutator  and assume that
$$
  \int_{V'} \langle A^{-1} u_1,u_1\rangle_{h_1} dV_{\omega} < \infty
$$
Then there exist an $(n,q-1)$-form $v_1$ with
$$
  \|v_1\|^2 =   \int_{V'} \langle v_1,v_1\rangle_{h_1} dV_{\omega} \leq \int_{V'} \langle A^{-1} u_1,u_1\rangle_{h_1} dV_{\omega}.
$$
\end{pos}
\begin{pos} If $u$ is smooth and $v_1$ is the minimal norm solution then  $v_1$  and the associated $(p,q-1)$-form $v$ are smooth.
The $L^2$-norms of $v_1$ on the sets
$V'_{\delta}:= \lbrace z \in V', d(z,\pi^{-1}(S))> \delta\rbrace$  with respect to  $h_Z$  and $h_E$  grow at most polynomially with respect to $\delta$ as $\delta \ra 0$ and the same holds  for  the $L^2$-norms of the corresponding $(p,q-1)$-form $v$.
\end{pos}
If $q = 1$ we can use Lemma 4.5 in \cite{pre2} which is an adaptation  of  Lemma 3.2 in \cite{fl} to get sup-norm estimates from the Bochner-Martinelli-Koppelman formula if we take a slightly smaller Stein neighbourgood $V'' \subset V',$ conic along $a(N(g))$ (Fig. \ref{okoliceb})
\begin{pos} If $q = 0$ and the initial form $u$ is smooth on a neighbourhood of $a(\ol{U}\setminus S),$ the form $v$ has at most polynomial poles on
  $\pi^{-1}(S).$
\end{pos}

\noi{\bf Notation.}
The notation from Theorem \ref{main thm1} is fixed throughout the paper. Let  $h_Z$ be a  Hermitian metric defined on the manifold $Z$ and let $\sigma: E \ra Z$ be a holomorphic vector bundle of rank $r$  equipped with a Hermitian metric $h_E.$
 The local coordinate system in a neighbourhood $V_{z_0} \subset Z$ of  a point $z_0 \in a(U)$ is $(z,w),$ where $z$ denotes the horizontal and $w$ the vertical (or fibre) direction and $z_0 =(0,0).$ More precisely, every point in $a(U)$ has $w = 0$ and points in the same fibre have the same first coordinate. If the point  $z_0$ is in $a(S)$ we write the $z$-coordinate as $z = (z_1,z_2),$ where $a(S)\cap V_{z_0} = \lbrace z_2 = 0, w = 0 \rbrace \cap V_{z_0}.$ The manifold $Z$ is $n$-dimensional and the dimension of the fibres $Z_{z_0}$ is constant, $r_0 = \dim Z_{z_0}.$ The notation $\zeta_1, \ldots, \zeta_n$ is sometimes used for local coordinates in $Z.$


\section{Construction of the K\"{a}hler metric $\omega$}\label{kahform}

The K\"{a}hler metric $\omega$  will be  obtained from the K\"{a}hler potential $\Phi = \vph_0 + \vph_1$ and is similar to the one constructed in  subsection 2.1 in \cite{pre2}. The only difference is that we choose a specific  \psh function $\vph_0$ instead of the given $\vph$ in  order to be able to study the curvature properties of $\omega.$ The construction is explained below.

 Since the Remmert reduction $p : X \ra \hat X$ of $X$ is a Stein space with finitely many isolated singular points it has a proper holomorphic embedding $\hat{f} :\hat{X} \ra \Cc^M$ for some large $M.$ Consequently, the holomorphic functions $ \hat{f}_{1} \circ p,\ldots,  \hat{f}_{M}\circ p$ generate the cotangent space $T^*{X \setminus S}$ and the function $\hat{\vph}_0:= \sum |\hat{f}_{i} \circ p|^2$ is a plurisubharmonic exhaustion function of $X$, strictly \psh on $X \setminus S.$ The functions $f_i:= \hat{f}_i \circ p \circ \pi $ are defined on $Z$ and they generate the horizontal cotangent space on $a(X\setminus S).$ Define $\vph_0 := \hat{\vph}_0 \circ \pi = \sum |f_i|^2$ and denote by $k_0$ the maximal order of degeneracy of  $i\partial \dbar\vph_0$ at $\pi^{-1}(S).$

For given $l > 2$ the construction in Subsection 2.1 in \cite{pre2} with the function $\vph$ replaced by $\vph_0$  yields almost holomorphic functions $f_{M+1},\ldots f_N,$ holomorphic to a  degree $l $  with zeroes of order at least $k$ on $\pi^{-1}(S).$  To be precise,  for every sufficiently large $k$ by Theorem A there exist sections $F_{M+1},\ldots, F_N  \in \Gamma(U',{\mathcal J}(a(S))^{k} ({\mathcal{J}}(a(U'))/{\mathcal J}^{l+1}(a(U'))))$, $U \Subset  U'$ which locally generate the sheaf on $\ol{U}.$ We further assume that $k > k_0.$
The functions $f_{M+1},\ldots f_N,$ are obtained by patching together particular local lifts of these sections using the partition of unity  $\lbrace \chi_j, U_j \rbrace$ which we (can) choose to depend on the horizontal variables only and so  $f_{M+1},\ldots f_N,$ are  holomorphic in  vertical directions. In local coordinates  $(z,w)$ near $a(S)$ we thus have
$$
  f_i(z,w) = \sum_{|\va| = k, 0 < |\vb| \leq l} c_{\va \vb} z_2^{\va} w^{\vb} +  \sum_{i,j,l} \chi_j(z) f_{i j l}(z,w)
$$
with $f_{i j l} \in {\mathcal O}(\|z_2\|^{k}\|w\|^{l+1})$ holomorphic on open sets $U_j.$
The functions $f_i$ satisfy the  estimate $\dbar f_i(z,w) \approx \|z_2\|^k \|w\|^{l+1}.$
 Moreover,  we can express  $z_2^{\va} w_j  = \sum g_{\va i j}(z) f_{i}(z,w) + {\mathcal O}(\|z_2\|^{k}\|w\|^{l+1})$ with $g_{\va i j}$ holomorphic and from this we infer that
 $\partial_{w_j}f_i,\, 1 \leq j \leq r_0, M+1 \leq i \leq N$
 generate the vertical cotangent bundle on a neighbourhood $V_T$ of $a(\ol{U})$ in $Z$ except on $\pi^{-1}(S).$
 Consequently, the matrix corresponding to $\partial_w \dbar_w \sum |f_i|^2$ is of the form $\|z_2\|^{2k} G,$ with $G$  invertible.
Define
$$
  \Phi:= \sum f_i \ol{f_i} \mbox{   and   } \omega:=i\partial \dbar \Phi.
$$
We claim that the function $\Phi$ is  a K\"{a}hler potential on a conic neighbourhood of $a(U\setminus S)$ and
$\omega$ a  K\"{a}hler metric.

Write  local coordinates  as $(\zeta_1,\ldots\zeta_n) = (z,w)$ and  represent the Levi form
$$
  i\partial \dbar \Phi = i\sum h_{jk} d\zeta_j \wedge d\ol{\zeta_k}
$$ by a matrix $H = \lbrace h_{jk}\rbrace.$
In local coordinates $(z,w)$ the nonnegative part of $\omega,$  $\omega_+ = i\sum \partial f_i \wedge \ol{\partial f_i},$ represented in the matrix form as $H_+,$ can be estimated from below by
\begin{equation}\label{h+}
  H_+(z,w)  \gtrsim   \begin{bmatrix}
    \|z_2\|^{2k_0}  + \|w\|^2 \|z_2\|^{2k-2}& \|w\|\|z_2\|^{2k-1}\\
     \|w\|\|z_2\|^{2k-1} & \|z_2\|^{2k}
\end{bmatrix},
\end{equation}
where we have estimated the decay of  $\vph_0$ by $\|z_2\|^{2k_0}$ from below.
The possibly negative part $\omega_{-}=i\sum \partial \dbar f_i \ol{f_i} + f_i \partial \dbar \ol{f_i} + \partial \ol{ f_i} \wedge \dbar{ f_i}$ degenerates at least as
$\|w\|^{l}\|z_2\|^{k-1}.$  It is clear that for $\|z_2\| > \delta $ and small $\|w\|$ or  $\|z_2\| \leq \delta $ and $\|w\| \leq \|z_2\|^2$ the matrix $H$ is strictly positive definite and thus $\omega$ is  a K\"{a}hler metric on a neighbourhood of $a(U \setminus S),$ conic along $a(S).$
Let   $(z,w)$ be local coordinates near $a(S)$ and define $H_0(z):= H(z,0)=H_+(z,0),$  $H_1 = H - H_0.$ It follows that $H_1 = \m O(\|w\|\|z_2\|^{2k-1})$ and
that $H_0$ decreases polynomially (in some directions) as we approach $\pi^{-1}(S)$ and  its degeneracy is bounded from below by $\|z_2\|^{2k},$
$$H_0(z) \gtrsim  \begin{bmatrix}
     \|z_2\|^{2k_0} & 0\\
     0 & \|z_2\|^{2k}
\end{bmatrix}.
$$
Notice that $H_0$ is strictly positive on a neighbourhood  of $a(\ol{U})$ except on $\pi^{-1}(S)$ (and therefore invertible) and
$\|H_0^{-1}\|$ degenerates in the worst case as $\|z_2\|^{-\kappa}$ for some $\kappa \geq 0.$ Because $S$ is compact, there exists one $\kappa$
for all points in $a(S).$
Write  $H = H_0(I + H_0^{-1} H_1),$   $H^{-1} = (I + H_0^{-1} H_1)^{-1}H_0^{-1},$ then
$$
  \|H_0^{-1}H_1\|\approx\|z_2\|^{2k-1-\kappa}\|w\|
$$
in the worst case and this term is small, $\|H_0^{-1}H_1\| < \|z_2\|^{3\kappa + k_1 }$ on conic neighbourhoods of the form $\|w\| \leq \|z_2\|^{4\kappa  + k_1},$ and so
\begin{equation}\label{o0}
  H^{-1} = H_0^{-1} +  H_0^{-1} H_1\sum_0^{\infty}(H_0^{-1}H_1)^n H_0^{-1} = H_0^{-1} + N,\, \|N\| \leq \|z_2\|^{2\kappa + k_1}.
\end{equation}
Inside this cone the  the degeneracy of the inverse $H^{-1}$  is governed by $H_0^{-1}.$

\section{Theorems on curvatures}
\subsection{Basic theorems on curvatures}\label{osnovni}

Before proceeding to the proof we recall some formulae from Demailly's  Complex analytic and algebraic geometry \cite{dem}.

Let $(X,\omega)$ be a K\"{a}hler manifold, $E \ra X$ a rank $r$ vector bundle equipped with a Hermitian metric $h.$ The matrix $H$ that corresponds to $h$ in local coordinates is given by $ \langle u,v \rangle_h = \sum h_{\lambda \mu} u_{\lambda} \ol{v_\mu} = u^T H \ol{v}.$
Let $i\Theta(E)$ be the Chern curvature form of the metric  and $\Lambda$ the adjoint of the operator
$u \ra u \wedge \omega,$ defined on $(p,q)$-forms. Denote by $L^2_{p,q} (X,E)$ the space of $(p,q)$-forms with bounded
$L^2$-norms with respect to the  $h$ and let
$ A_{E,\omega} = [i\Theta(E), \Lambda]$ be the commutator.

In bidegree $(n,q)$  the positivity of  $A_{E,\omega}$ is equivalent to Nakano positivity of $E.$ Let $e_1,\ldots,e_r$ be a local frame of $E.$ If
the metric is locally represented by a matrix $H$ then
\begin{equation}\label{forma}
 i\Theta(E) = i\dbar(\ol{H}^{-1} \partial \ol{H})=i\sum c_{j k \lambda \mu} dz_j\wedge d \ol{z}_k \otimes e^*_{\lambda}\otimes e_{\mu},
\end{equation}
If  $e_1,\ldots,e_r$ is an orthonormal frame then  the Hermitian form $\theta_E$ defined on $TX \otimes E,$ which is associated to $i\Theta(E),$ takes the form
\begin{equation}\label{forma1}
  \theta_E = \sum  c_{j k \lambda \mu} (dz_j \otimes e^*_{\lambda}) \otimes \ol{(d{z}_k  \otimes e^*_{\mu})}.
\end{equation}
The curvature tensor (\ref{forma}) is {\it Griffiths positive} if the form (\ref{forma1}) is positive on  decomposable tensors $\tau = \xi \otimes v,$ $\xi \in TX,$ $v \in E,$
$\theta_E(\tau,\tau) = \sum c_{j k \lambda \mu} \xi_{j} \ol{\xi}_{k} v_{\lambda} \ol{v}_{\mu}$ and
{\it Nakano positive} if it is positive on $\tau = \sum \tau_{j \lambda} ({\partial}/{\partial z_j}) \otimes e_{\lambda},$
$\theta_E(\tau,\tau) = \sum c_{j k \lambda \mu} \tau_{j\lambda} \ol{\tau}_{k\mu}.$
In a nonotrhonormal frame we have
$\theta_E(\tau,\tau) = \sum c_{j k \lambda \mu} \tau_{j\lambda} \ol{\tau}_{k\nu} h_{\mu \nu}.$ Proposition VII-9.1,\cite{dem} states that
\begin{equation}\label{zveza}
\mbox{if  } \theta(E) >_{\rm Grif} 0 \mbox{ then }  r\Tr_E(\theta(E)) \otimes h - \theta(E) >_{\rm Nak} 0.
\end{equation}

The metric $h$ on $E$ induces the metric $h^s$ on $\Lambda^sE.$  Let $L$ be an $s$-tuple of (not necessarily ordered)  indices $L = (\lambda_1,\ldots,\lambda_s)$ and denote  $e_{L} := e_{\lambda_1}\wedge \ldots \wedge e_{\lambda_s}.$ If $\sigma$ is a permutation then $e_{\sigma(L)}=\sgn(\sigma)e_{L}.$  Let $L, M\in \lbrace (\lambda_1,\ldots \lambda_s), 1 \leq  \lambda_1 < \ldots < \lambda_s \leq r \rbrace=: {\mathcal L} $, $L =(\lambda_1,\ldots \lambda_s),$ $M =(\mu_1,\ldots,\mu_s).$
The coefficient $H_{LM} = \langle e_{L}, e_{M} \rangle_{h^s}$ in the matrix $H^s$ representing the induced metric $h^s$ is
$$
  H^s_{LM} =  \det H_{(\lambda_1,\ldots \lambda_s),(\mu_1,\ldots,\mu_s)},
$$
where $H_{(\lambda_1,\ldots \lambda_s),(\mu_1,\ldots,\mu_s)} $ is a submatrix of $H$ generated by rows $\lambda_1,\ldots \lambda_s$  and  columns $\mu_1,\ldots,\mu_s $ of the matrix $H.$
If  $e_1,\ldots e_r$ are orthonormal at $z$ so are their wedge products $\lbrace e_{L}, L \in {\mathcal L} \rbrace.$

The induced Chern curvature tensor on $\Lambda^s E, $  $i\Theta(\Lambda^s(E))  = \sum_{j,k} i\Theta(\Lambda^s(E))_{j k}dz_j \wedge d\ol{z}_k$
is defined by formula V-(4.5'), \cite{dem},
\begin{equation}\label{forma_s}
  i\Theta(\Lambda^s(E))_{jk}(e_{L}) = i\sum_{1 \leq l \leq s}  e_{\lambda_1} \wedge \ldots \wedge\Theta(E)_{jk} e_{\lambda_l}\wedge \ldots \wedge e_{\lambda_s}.
\end{equation}
 It is known that $E \geq_{\rm Nak} 0$ implies $\Lambda^s E \geq_{\rm Nak} 0.$   The following lemma gives an
explicit formula for the curvature  $i\Theta(\Lambda^s(E))$ in terms of the curvature $i\Theta(E)$ and shows that if
the associate Hermitian form $\theta(E)$ has at most polynomial poles on $\pi^{-1}(S),$ so does $\theta(\Lambda^s(E)).$

\begin{lm}\label{rformula} If $i\Theta(E)$ is Nakano nonpositive (nonnegative) then $i\Theta(\Lambda^s(E)),$ $1 \leq s \leq r,$ is also Nakano nonpositive (nonnegative).
\end{lm}

\dok
By  formula (\ref{forma_s}) we have
$$
  i\Theta(\Lambda^s(E))(e_{L})_{jk} = i\sum_{l, \mu} (-1)^{l-1}c_{jk \lambda_l \mu}  e_{\mu} \wedge e_{L_l'},
$$
where $L_l'$ is obtained from $L$ by removing the $l$-th index.
Let $L(\lambda,\mu)$ denote the (not ordered) multiindex  obtained by replacing the index $\lambda$ in the multiindex $L$  by  $\mu.$
We define that $e_{L(\lambda, \mu)} = 0$ if and only if $\lambda \notin L$ or $\mu \in L\setminus \lbrace \lambda \rbrace.$
If $\lambda_l = \lambda$ then $e_{\mu L_l'} = (-1)^{l-1}e_{L(\lambda,\mu)}$ and
\begin{eqnarray*}
  i\Theta(\Lambda^s(E)) &=& i\sum_{j,k,L,M} c^s_{jk L M} dz_j\wedge d{\ol{z}}_k \otimes e_{L}^* \otimes e_M \\
    &=& i\sum_{j,k,L} \; \sum_{\lambda \in L,\mu} c_{jk\lambda\mu} dz_j\wedge d{\ol{z}}_k \otimes e_{L}^* \otimes e_{ L(\lambda,\mu)}\\
 &=& i\sum_{\substack{j,k\\
               |L'| = s-1}} \sum_{\lambda,\mu \not\in L'} c_{jk\lambda\mu} dz_j\wedge d{\ol{z}}_k \otimes e_{\lambda L'}^* \otimes e_{\mu L'}.
 \end{eqnarray*}
Here the bijection between the sets $\lbrace L : |L| = s \rbrace$
and $\lbrace \lambda L' : |L'| = s-1, \lambda \notin L'\rbrace$ is used, where we compare multiindices as sets.
Let $\tau = \sum_{j,L}\tau_{j,L} (\partial/\partial z_j) \otimes e_{L}$ with  additional properties $\tau_{j,\sigma(L)} := \sgn(\sigma)\tau_{j,L}$ for any permutation $\sigma$ and $ \tau_{j,L} = 0$ if there are at least two equal indices in $L.$
The bundle $\Lambda^s(E)$ is Nakano positive if the bilinear form
\begin{eqnarray*}
  \theta_{\Lambda^s(E)}(\tau,\tau) &=& \sum c^s_{jk L M} \tau_{j,L} \ol{\tau}_{k,S} \langle e_{M},e_{S}\rangle_{h^s}\\
     &=& \sum_{\substack{
         j,k,|S|=s\\  |L'|=s-1 }} \sum_{\lambda,\mu \not\in L'} c_{jk\lambda\mu} \tau_{j,\lambda L'} \ol{\tau}_{k,S}\langle e_{\mu L'} ,e_{S}\rangle_{h^s}
\end{eqnarray*}
is positive.
Assume that the local frame $e_1, \ldots, e_r$ is orthonormal. Then $\lbrace e_{L}, L \in {\mathcal L} \rbrace$ are orthonormal and
\begin{equation}\label{e1}
  \theta_{\Lambda^s(E)}(\tau,\tau) =  \sum_{j,k,L'} \sum_{\lambda,\mu} c_{jk\lambda\mu}  \tau_{j,\lambda L'}\ol{\tau}_{k,\mu L'}= \sum_{L'} \theta_E(\tau_{L'},\tau_{L'}),
\end{equation}
 where the form $ \tau_{L'} $ is defined by
$\tau_{L'}=\sum_{j,\lambda} \tau_{j,\lambda L'}(\partial/\partial z_j) \otimes e_{\lambda} $
for any multiindex $L'$ of length $s-1.$
Hence  if $\theta_E$ is Nakano nonpositive (nonnegative), so is $\theta_{\Lambda^sE}.$
\qed


\subsection{Almost nonpositivity  of $\omega$}

In this section we study properties of the form $\omega$ constructed in Sect. \ref{kahform}. Let $V \subset Z$ be a neighbourhood of $a(U \setminus S),$ conic along $a(S).$
The metric $\omega$ on a vector bundle $E \ra V$ is {\it almost Nakano nonpositive} if the curvarure tensor $i\Theta(E)$ has a decomposition
$i\Theta(E) = i\Theta_0(E) + i\Theta_1(E),$ where $i\Theta_0(E)$ is nonpositive and $i\Theta_1(E)$ is locally of the form
$i\Theta_1(E)(z,w) = {\mathcal O}(\|w\|^l\|z_2\|^k)$ near points in $a(S)$ and $i\Theta_1(E)(z,w) = {\mathcal O}(\|w\|^l)$ near points in $a(U \setminus S)$, for some $l \in \Nn,$ $k \in\Zz$. \\

\noi The main theorem in this subsection is

\begin{izr}[Almost Nakano nonpositive K\"{a}hler metric]\label{ANK}  Let $Z,X,S,a,U$  and $\omega$  be as in Theorem \ref{main thm1}.
 There exist a neighbourhood $V$ of $a(\ol{U} \setminus S)$ conic along $a(S)$ such that the metric $\omega$ on $TZ|_{V}$
is almost Nakano nonpositive.
\end{izr}

\begin{pos}\label{negativni} Let $h_{\omega}$ be the metric on $TZ|_V$ induced by $\omega.$ Then $h_{\omega} e^{ \Phi}$ is Nakano negative on
a smaller  neighbourhood - which we again denote by  $V$ - of $a(U \setminus S),$ conic along $a(S).$
\end{pos}

\noi{\it Proof of Theorem \ref{ANK}}.
Write $f = (f_1,\ldots f_N)^T$  and  $\sum|f_i|^2 = f^T\ol{f}$ and let
 $H$ denote the matrix corresponding to the metric.
If $D$ denotes the holomorphic and $\ol{D}$ the antiholomorphic derivative with respect to
$(z,w),$
$$
   Df =  \begin{bmatrix}
    f_{1,z_1} & f_{1, z_2} & f_{1, z_3} & \dots  & f_{1,w_{r_0}} \\
    f_{2,z_1} & f_{2,z_2} & f_{2,z_3} & \dots  & f_{2w_{r_0}} \\
    \vdots & \vdots & \vdots & \ddots & \vdots \\
    f_{N,z_1} & f_{N,z_2} & f_{N,z_3} & \dots  & f_{N,w_{r_0}}
\end{bmatrix}
$$
then
$$
    H = D \ol{D} (f^T \ol{f}) =  (Df)^T\ol{Df} + D\bar{f}^T\ol{D}{f}+ f^T L\ol{f} +  Lf^T \ol{f}
$$
The Levi form $Lf$ of the vector is calculated as $Lf = (Lf_1,\ldots, Lf_N).$ 
  With the notation defined prior to (\ref{h+}) we have  $H_+ =  (Df)^T\ol{Df}$ and $H_- = H - H_+ = D\bar{f}^T\ol{D}{f}+ f^T L\ol{f} +  Lf^T \ol{f}. $
Since $f_i$ are holomorphic or almost holomorphic to the degree $l$, we have  $\dbar f = {\mathcal O}(\|w\|^{l+1}),$ terms $ \dbar Df, Lf,  H_-$ are of the form ${\mathcal O}(\|w\|^l)$  and because of the holomorphicity of $f_i$ in $w$-directions the terms
$\partial H_-, \dbar H_-,\partial\dbar Df,\partial \dbar \ol{H}_-$ are of the form ${\mathcal O}(\|w\|^{l-1}).$ The term $H_+$ gives
\begin{eqnarray*}
  \partial \dbar \ol{H}_+ &=& \partial \dbar ((Df)^*{Df}) = -\dbar \partial ((Df)^* {Df})= -\dbar (\partial (Df)^* {Df} + (Df)^*(\partial{Df})) \\
     &=& -(\partial {Df})^*\wedge(\partial {Df}) + \partial (Df)^*\wedge \dbar Df  + (\partial \dbar (Df)^*){Df} + (Df)^* \partial \dbar {Df}
\end{eqnarray*}
and so
\begin{equation}\label{o1}
  \partial  \dbar \ol{H} = -(\partial {Df})^*\wedge(\partial {Df})+ {\mathcal O}(\|w\|^{l-1}).\
\end{equation}
By the above estimates we conclude that
\begin{equation}
  \partial \ol{H} =\partial (\ol{H}_+ + \ol{H}_-)=(Df)^*(\partial{Df})  +  {\mathcal O}(\|w\|^{l-1}) \label{o2}
\end{equation}
and similarly
\begin{equation}
  \dbar\ol{H}  =(\partial Df)^* {Df} +  {\mathcal O}(\|w\|^{l-1}). \label{o3}
\end{equation}
Since the metric  is K\"{a}hler, we may assume that the coordinates $\zeta_1,\ldots,\zeta_n$ near
the point $z_0 \in a(U\setminus S)$ are such that $H(z_0) = I$ to the second order. The curvature form equals
\begin{eqnarray*}
   i\Theta(TZ) &=& i\dbar (\ol{H}^{-1}\partial\ol{H})  \\
      &=& -i\ol{H}^{-1}\dbar \ol{H} \ol{H}^{-1}\wedge \partial \ol{H} + i\ol{H}^{-1}\dbar\partial \ol{H}\\
      &=& -i\partial \dbar\, \ol{H} = i(\partial {Df})^*\wedge(\partial {Df}) + {\mathcal O}(\|w\|^{l-1})
\end{eqnarray*}
by the above assumptions. We claim that  $i(\partial {Df})^*\wedge(\partial {Df})$ is Nakano nonpositive. Denote the dual tangent vectors by
 $e_{\lambda} := \partial/\partial \zeta_{\lambda}$ and
let
$$
  i\Theta_0 = i(\partial {Df})^*\wedge(\partial {Df}) = i\sum_{j,k,\lambda, \mu}c_{jk\lambda\mu}d\zeta_j\wedge  d\ol{\zeta}_k \otimes e^*_{\lambda}\otimes e_{\mu}.
$$
Then we have
$$
   c_{jk\lambda \mu} = -\sum_{i = 1}^N \frac{\partial^2 f_i}{\partial \zeta_j \partial \zeta_{\lambda} }\ol{\frac{\partial^2 f_i}{\partial \zeta_k \partial \zeta_{\mu} }}
$$
and so
$$
  \theta_0(\tau,\tau) = - \sum_{i}\left|\sum_{ j \lambda} \frac{\partial^2 f_i}{\partial \zeta_j \partial \zeta_{\lambda} }\tau_{j\lambda}\right|^2 \leq 0.
$$

\begin{op} If the functions $f_i$ were holomorphic, then the Nakano nonpositivity of $i\Theta(VT)$ on $V$ could  be inferred from the fact that  the metric on  $TZ|_{V}$ is the metric induced on the subbundle $F\leq  V \times \Cc^N$ by the standard metric on $\Cc^N$  via the holomorphic vector bundle isomorphism $TZ_V \ra F,$ $TZ \ni v \mapsto (df_1(v),\ldots, df_N(v)) \in \Cc^N$ and it is known that the Nakano curvature decreases in subbundles by VII-(6.10), \cite{dem}. Because $f_i$ are almost holomorphic, we get an `error' term, denoted by $i\Theta_1$ in the sequel and we will show that it decreases arbitrarily fast on conic neighbourhoods.  On the section $a(U \setminus S)$ the `holomorphic' part of the curvature tensor is exactly $i \Theta_0.$
\end{op}

The curvature form $i\Theta(TZ)$ is then almost Nakano nonpositive
on an open neighbourhood of $a(U \setminus S).$ If we want to prove that the neighbourhood is conic we have to show that the
part of the curvature which contains $w$-variables  decreases polynomially sufficiently fast  in some conic neighbourhood. Let $(z,w)$ be local coordinates at $(z,0) \in a(S).$  We can not  assume that $H=I$ to the second order at $(z,0)$ because $H(z,0)$ is degenerate on $a(S).$

We have to split the curvature tensor into the part depending only on $z$-variables - we have just proved that it is nonpositive - and the rest, which we want to be small on conic neighbourhoods.
By estimates (\ref{o0}), (\ref{o1}), (\ref{o2}), (\ref{o3}) we have for $\|w\| \leq \|z_2\|^{4\kappa  + k_1}$
\begin{eqnarray*}
  \ol{H}^{-1}\dbar \partial \ol{H}&=& \ol{H}^{-1}((\partial {Df})^*\wedge(\partial {Df})+ {\mathcal O}(\|w\|^{l-1}) )\\
                       &=&\ol{H}_0^{-1}(\partial {Df})^*\wedge(\partial {Df}) + \ol{N} (\partial {Df})^*\wedge(\partial {Df})  +
                        {\mathcal O}(\|w\|^{l-1}\|z_2\|^{-\kappa})\\
                       &=& \ol{H}_0^{-1}(\partial {Df})^*\wedge(\partial {Df}) + \mathcal O(\|z_2\|^{2\kappa + k_1}) + {\mathcal O}(\|w\|^{l-1}\|z_2\|^{-\kappa}),
\end{eqnarray*}
\begin{eqnarray*}
   &&\ol{H}^{-1}\dbar \ol{H} \ol{H}^{-1}\wedge\partial \ol{H} =  (\ol{H}_0^{-1} + \ol{N})((\partial Df)^* {Df} + {\mathcal O}(\|w\|^{l-1}))
    (\ol{H}_0^{-1} + \ol{N})\\
    && \qquad \qquad \qquad \qquad\qquad \wedge((Df)^*(\partial{Df}) +{\mathcal O}(\|w\|^{l-1}))\\
   &&  \qquad =\ol{H}^{-1}_0 (\partial Df)^* {Df} \ol{H}^{-1}_0 \wedge(Df)^*(\partial{Df})
    +\ol{H}^{-1}_0 (\partial Df)^* {Df} \ol{N} \wedge(Df)^*(\partial{Df})\\
    &&\qquad \qquad +\ol{N} (\partial Df)^* {Df} \ol{H}^{-1}_0 \wedge(Df)^*(\partial{Df})
    +\ol{N}(\partial Df)^* {Df} \ol{N} \wedge(Df)^*(\partial{Df})\\
    &&\qquad \qquad + \mathcal O(\|z_2\|^{-2\kappa}) \mathcal O(\|w\|^{l-1}) \\
    && \qquad = \ol{H}^{-1}_0 (\partial Df)^* {Df} \ol{H}^{-1}_0\wedge (Df)^*(\partial{Df}) + \mathcal O(\|z_2\|^{2\kappa + k_1}) \mathcal O(\|z_2\|^{-\kappa})  \\
    &&\qquad \qquad + \mathcal O(\|z_2\|^{2(2\kappa + k_1)})  + \mathcal O(\|z_2\|^{-2\kappa}) \mathcal O(\|w\|^{l-1}).
\end{eqnarray*}
Let $\Theta_0(TZ) := \ol{H}_0^{-1}(\partial {Df})^*\wedge(\partial {Df}) - \ol{H}^{-1}_0 (\partial Df)^* {Df} \ol{H}^{-1}_0\wedge (Df)^*(\partial{Df}). $
By Taylor series expansion we see that
\begin{eqnarray*}
   i\Theta_0(TZ)(z,w) &=& i\ol{H}_0^{-1}(\partial {Df}(z,0) + \m O(\|w\|))^*\wedge(\partial {Df}(z,0) + \m O(\|w\|))\\
                     &&\qquad- i\ol{H}^{-1}_0 (\partial Df(z,0) + \m O(\|w\|))^* (Df(z,0) + \m O(\|w\|))\\
                     && \qquad \cdot \ol{H}^{-1}_0 \wedge (Df(z,0)+ \m O(\|w\|))^*(\partial{Df}(z,0) + \m O(\|w\|))\\
                     &=& i\Theta_0(TZ)(z,0) + \m O(\|w\|)\m O(\|z_2\|^{-\kappa}  + \|z_2\|^{-2\kappa}).
\end{eqnarray*}
Let $i\Theta_1(TZ)(z,w) = i\Theta(TZ)(z,w) - i\Theta_0(TZ)(z,0).$
The form $i\Theta_0(z,0)$ is nonpositive on a neighbourhood of $a(U\setminus S),$ conic along $a(S).$  The `error term' $i\Theta_1(TZ)$ is in the worst case
\begin{eqnarray*}
     &&\m O(\|w\|)\m O(\|z_2\|^{-\kappa}  + \|z_2\|^{-2\kappa}) + \mathcal O(\|z_2\|^{2\kappa + k_1}) + {\mathcal O}(\|w\|^{l-1}\|z_2\|^{-\kappa}) \\
    && + \mathcal O(\|z_2\|^{\kappa + k_1})
     + \mathcal O(\|z_2\|^{2(2\kappa + k_1)})  + \mathcal O(\|z_2\|^{-2\kappa}) \mathcal O(\|w\|^{l-1})
\end{eqnarray*}
 and it decreases at least as $\|z_2\|^{k_1}$  on conic neighbourhoods $\|w\| \leq \|z_2\|^{4\kappa  + k_1}.$\qed\\

\noi{\it Proof of Corollary \ref{negativni}}.
The part $H_0$ of the Levi form of $\Phi$ near $(z,0) \in a(S)$ that dominates in the matrix $H$ on a cone is bounded from below by $\|z_2\|^{2k}I.$
The potentially positive part of the Nakano curvature, $i\Theta_1,$ is of the form
$$
   {\mathcal O}(\|z_2\|^{k_1})
$$
on $\|w\| \leq \|z_2\|^{4\kappa  + k_1}$ and can be compensated by $-\partial\dbar \Phi/2$ for $k_1 > 2k$ thus making
the curvature tensor
$$
  (i\Theta_0(z,0) - i\partial \dbar \Phi(z,w)) + (i \Theta_1(z,w)  - i\partial \dbar \Phi(z,w)))
$$
strictly Nakano negative. \qed

\section{Proof of the main theorem}



By  theorem 1.1 in \cite{pre2} the bundle $E$ can be endowed with a Nakano positive Hermitian metric $h_0$ on a conic neighbourhood of $a(U \setminus S)$ with
polynomial poles on $\pi^{-1}(S).$
Let $\omega = i\partial\dbar \Phi$ be the given metric.
 In order to solve the $\dbar$-equation in  bidegree $(p,q)$
we have to show that the curvature tensor
$$
  i\Theta(E \otimes \Lambda^s TZ) + iL\Psi = i\Theta(E) + i\Theta(\Lambda^s TZ)  +  iL\Psi
$$
is positive (or at least nonnegative) for some strictly \psh weight $\Psi$ and $s = n-p.$

The K\"{a}hler metric $h$ induced by the K\"{a}hler form  $\omega$ has almost nonpositive Nakano curvature. Let $h^s$ be the metric on $\Lambda^sTZ$ induced by $h$ and  $h_1$ the metric induced by the form $\omega_1 =\omega e^{ \Phi}.$ The latter has strictly Nakano negative curvature tensor
$$
  i\Theta(TZ)_{\omega_1} = i\Theta(TZ)_{\omega} - i\partial\dbar \Phi
$$
by Corollary
\ref{negativni}.
If the original metric on $TZ$ is represented by $H$ then the new one is $H_1 = H e^{\Phi}$ and the induced metric $h_1^s$ on $\Lambda^sTZ$ is represented by the matrix
$$
  H_1^s = H^s e^{{s}\Phi }.
$$
Since the Chern curvature tensor of the Hermitian metric $\omega_1$ on $TZ$ is Nakano negative   the induced curvature tensor
on $\Lambda^sTZ_{\omega_1}$ is also Nakano negative by Lemma \ref{rformula} and equals
$$
  i\Theta(\Lambda^sTZ)_{h_1^s} =  i\Theta(\Lambda^sTZ)_{h^s} - i{s}\partial\dbar \Phi <_{\rm Nak} 0.
$$
Write $F:=\Lambda^sTZ$ and
let $\theta_F$ be the bilinear form  on $TZ \otimes F $  associated to $i\Theta(F).$
Since the rank of the bundle $F$ is $\binom{n}{s},$ formula (\ref{zveza}) gives
$$
  \theta_1 =-\binom{n}{s}\Tr_F(\theta_F)_{h_1^s}\otimes h_1^s + \theta_{F_{h_1^s}} >_{\rm Nak} 0.
$$
We observe that $\theta_1$ is the curvature form associated to the Chern curvature tensor of
the bundle $F \otimes( \det F^*)^{\binom{n}{s}}$ with the metric induced by $h_1.$
The induced metric on $\det F^*$ equals
\begin{equation}\label{dual}
  h_{\det F^* } = (\det H_1)^{-\binom{n-1}{s-1}} =(\det H)^{-\binom{n-1}{s-1}}e^{-\binom{n-1}{s-1}\Phi}
\end{equation}
by the identity
\begin{equation}\label{id}
  \det F = \det \Lambda^sTZ = (\det TZ)^{\binom{n-1}{s-1}}
\end{equation}
and so
$$
  i\Theta(\det F^*)_{h_1} = i\Theta(\det F^*)_{h} + \binom{n-1}{s-1}i \partial\dbar \Phi.
$$
Then
$$
  \theta_1 =\theta((\det F^*)^{\binom{n}{s}})_{h}\otimes h_1^s + e^{s \Phi} \theta(F)_{h} + \left({\binom{n}{s}}{\binom{n-1}{s-1}} - s\right)\partial\dbar\Phi \otimes h_1^s >_{\rm Nak} 0
$$
or
\begin{equation}\label{nak_pos}
  \theta =\theta((\det F^*)^{\binom{n}{s}})_{h}\otimes h^s +  \theta(F)_{h} + \left({\binom{n}{s}}{\binom{n-1}{s-1}} - s\right)\partial\dbar\Phi \otimes h^s >_{\rm Nak} 0.
\end{equation}
To complete the proof we  view this expression as part of the curvature of the metric $e^{-\Psi} h^s $ on $\Lambda^sTZ.$
Let $\Phi_1:=({\binom{n}{s}}{\binom{n-1}{s-1}} - s)\Phi.$ Then the weight $e^{-\Phi_1}$ gives the last term of $\theta.$
Observe that
$$
  i\Theta(\det F^*)_{h} = i\dbar\partial \log(\det H)^{-\binom{n-1}{s-1}} = i\partial\dbar \log(\det H)^{\binom{n-1}{s-1}}.
$$
Write
$$
  F = F \otimes( \det F^*)^{\binom{n}{s}} \otimes (\det F)^{\binom{n}{s}}=F \otimes( \det F^*)^{\binom{n}{s}} \otimes (\det TZ) ^{\binom{n-1}{s-1}\binom{n}{s}}
$$
by invoking (\ref{id}) and define the Nakano positive metric on $(\det TZ)^{\binom{n-1}{s-1}\binom{n}{s}}$ in the following way.
Let $v_i$ be smooth sections of $\det TZ,$ given by Proposition 2.1 in \cite{pre2}, which are holomorphic to the degree $l_2>2$ with zeroes of order $k_2$ on $\pi^{-1}(S)$ and such that they  generate the bundle $\det TZ$ on a neighbourhood of $a(U\setminus S),$ conic along $a(S).$ Let $v$ be a local holomorphic section, $v_i = \va_i v$  and define
$$
  \Phi_2:= \log \sum \langle v_i,v_i \rangle_{H} = \log \langle v,v \rangle_{H} + \log \sum |\va_i|^2.
$$
Then the metric $e^{-\Psi}h^s$ for $\Psi = \Phi + \Phi_1 + \binom{n}{s}\binom{n-1}{s-1}\Phi_2$ has
\begin{eqnarray*}
&& i\binom{n}{s}\binom{n-1}{s-1}\partial\dbar\log\det H+ i\Theta(\Lambda^sTZ)_{\omega}+ i\left( {\binom{n}{s}}{\binom{n-1}{s-1}} - s \right) \partial\dbar\Phi\\
&& \;\;\;
              + \,i\partial\dbar \Phi + i\binom{n}{s}\binom{n-1}{s-1}\partial\dbar\log \sum |\va_i|^2
\end{eqnarray*}
as a curvature tensor. The first three terms give   a Nakano positive curvature by (\ref{nak_pos}) and the last is also
Nakano positive in a suitable conic neighbourhood because the negative part of $\partial \dbar \log \sum |\va_i|^2$ is of the form
$$
  C_1\frac{\|w\|^{l_2}}{\|z_2\|} + C_2\|w\|^{2l_2} + C_3 \|w\|^{l_2} + C_4{\|w\|^{l_2 - 1}}
$$
and its modulus decreases at least as $\|z_2\|^{k_1}$ on  conic neighbourhoods $\|w\| \leq \|z_2\|^{k_1+2}$ (see \cite{pre2}, p.14 for details). Because $i\partial \dbar\Phi$ is strictly plurisubharmonic with the rate of degeneracy  at most $\|z_2\|^{2k}$  (independent of the shape of the cone if it is sharp enough) it compensates the negativity of $\partial \dbar \log \sum|\va_i|^2$ provided $k_1 > 2k.$\qed

 \end{document}